\documentstyle[12pt]{amsart}
\def\ds{\displaystyle}
\title{A generalization of an inequality from IMO 2005}
\author{Nikolai Nikolov}
\address{Institute of Mathematics and Informatics\\ Bulgarian Academy of
Sciences\\ Acad. G. Bonchev Str., Block 8, 1113 Sofia, Bulgaria}
\email{nik@@math.bas.bg}
\begin{document}
\maketitle

The present paper was inspired by the third problem from the IMO
2005. A special award was given to Yurie Boreiko from Moldova for
his solution to this problem. It was the first such award in the
last ten years. Here is the problem.

\smallskip

{\bf Problem.} {\it Let $x$, $y$ and $z$ be positive real numbers
such that $xyz\geq 1$. Prove that $$\frac{x^5-x^2}{x^5+y^2+z^2} +
\frac{y^5-y^2}{y^5+z^2+x^2} + \frac{z^5-z^2}{z^5+x^2+y^2} \ge0. $$
}

The main objective of this paper is to prove the following more
general inequality:
\smallskip

{\bf Proposition 1.} {\it Let $x_1,x_2,\dots,x_n$ be positive real
numbers such that $\ds\prod_{i=1}^n x_i\ge 1.$ Then, for all
$\alpha\ge 1,$
$$\sum_{i=1}^n\frac{x_i^\alpha-x_i}{x_1+\dots+x_{i-1}
+x_i^\alpha+x_{i+1}+\dots+x_n}\ge 0.\eqno{(1)}$$}

{\bf Remark.} We get the mentioned IMO problem by choosing $n=3$
and $\ds\alpha=\frac52,$ and applying the result to the numbers
$x^2,$ $y^2$ and $z^2.$

The statement is trivial for $n=1;$ thus, we assume that $n\ge 2$.
We will consider two cases, depending on $\alpha.$ In the first
case we assume that $\ds\alpha\le2+\frac{1}{n-1}$. We can then
prove a stronger inequality in which all denominators are equal to
the sum of the numbers. This was the approach used by Yurie
Boreiko. (This idea was also suggested by the Armenian deputy
leader Nairi Sedrakyan.) The stronger inequality fails for
$\ds\alpha>2+\frac{1}{n-1}.$ In this case, we will estimate the
terms of the sum from below by suitable real numbers which sum to
zero and have equal denominators. Thus, we use two different ideas
for the two cases. The third problem from IMO 2205 is on the
boundary of the toe cases in the proof of Proposition 1. It would
be interesting to find a unified and concise approach for the
proof of Proposition 1.

{\it Proof of Proposition 1:}
\smallskip

{\bf Case 1.} $\ds1\le\alpha\le2+\frac{1}{n-1}.$

We have $$\frac{x_1^\alpha-x_1}{x_1^\alpha+\sum_{i=2}^n x_i}
=\frac{ x_1-\frac{1}{x_1^{\alpha-2}}}{x_1+\frac{\sum_{i=2}^n
x_i}{x_1^{\alpha-1}}} \ge\frac{
x_1-\frac{1}{x_1^{\alpha-2}}}{\sum_{i=1}^n x_i}\eqno{(2)}.$$ (To
prove the inequality in (2), consider the cases $x_1\ge1$ and
$x_1\le1$ separately). Using (2) and the analogous inequalities
for $x_2,\dots,x_n,$ we get (1) from the inequality
$\ds\sum_{i=1}^n\left(x_i-\frac{1}{x_i^{\alpha-2}}\right)\ge 0,$
that is, $$\sum_{i=1}^nx_i\ge\sum_{i=1}^nx_i^\beta,\eqno{(3)}$$
where $\ds\beta=2-\alpha\in\left[-\frac{1}{n-1},1\right].$ Now, to
prove (3) we will consider two subcases.
\smallskip

{\it Subcase 1.1.} $\beta\in[0,1].$

Then $x^\beta$ (for $x\ge 0$) is concave, and Jensen's Inequality
implies $$\left(\frac{\sum_{i=1}^nx_i}{n}\right)^\beta\ge
\frac{\sum_{i=1}^nx_i^\beta}{n}.\eqno{(4)}$$ On the other hand,
$\ds\prod_{i=1}^nx_i\ge 1,$ and therefore,
$\frac{\sum_{i=1}^nx_i}{n}\ge 1$ by the AM-GM Inequality. Since
$\beta\le1,$ we get $$\frac{\sum_{i=1}^n x_i}{n}\ge
\left(\frac{\sum_{i=1}^nx_i}{n}\right)^\beta.\eqno{(5)}$$ Now (3)
follows from (4) and (5).
\smallskip

{\it Subcase 1.2.} $\ds\beta\in\left[\frac{1}{1-n},0\right].$

Then $$x_1^\beta\le\ds\prod_{i=2}^n x_i^{-\beta}\le\frac{
\sum_{i=2}^n x_i^{\beta(1-n)}}{n-1}$$ by the AM-GM inequality.
Adding to this inequality the analogous inequalities for
$x_2,\dots,x_n$ yields
$$\sum_{i=1}^nx_i^\beta\le\sum_{i=1}^nx_i^{\beta(1-n)}.\eqno{(6)}$$
But $0\le\beta(1-n)\le 1.$ Thus
$$\sum_{i=1}^nx_i^{\beta(1-n)}\le\sum_{i=1}^nx_i\eqno{(7)}$$ by
Subcase 1.1. Now (3) follows from (6) and (7).
\smallskip

{\bf Case 2.} $\ds\alpha\ge2+\frac{1}{n-1}.$

It suffices to show that
$$\frac{x_1^\alpha-x_1}{x_1^\alpha+\sum_{i=2}^nx_i}\ge
\frac{nx_1^\gamma-\sum_{i=1}^n x_i^\gamma}{(n-1)\sum_{i=1}^n
x_i^\gamma}\eqno{(8)}$$ for some $\gamma,$ and then to add to (8)
the analogous inequalities for $x_2,\dots,x_n.$

Subtracting 1 from both sides in (8) gives $$-\frac{\sum_{i=1}^n
x_i}{x_1^\alpha+\sum_{i=2}^n x_i}\ge-\frac{n\sum_{i=2}^n
x_i^\gamma}{(n-1)\sum_{i=1}^n x_i^\gamma}$$ or, equivalently,
$$n\frac{x_1^\alpha+\sum_{i=2}^n x_i} {\sum_{i=1}^n
x_i}\ge\frac{(n-1)\sum_{i=1}^n x_i^\gamma}{\sum_{i=2}^n
x_i^\gamma}.$$ Subtracting $n$ from both sides yields
$$\frac{nx_1(x_1^{\alpha-1}-1)}{\sum_{i=1}^n
x_i}\ge\frac{(n-1)x_1^\gamma}{\sum_{i=2}^n x_i^\gamma}-1.$$ Since
$\ds\prod_{i=1}^nx_i\ge 1,$ we have $\ds\prod_{i=1}^n
x_i^{\frac{\alpha-1}{n}}\ge 1,$ and the above inequality will
follow from the homogeneous inequality $$\frac{nx_1}{\sum_{i=1}^n
x_i}\left(\frac{x_1^{\alpha-1}}{\ds\prod_{i=1}^n
x_i^{\frac{\alpha-1}{n}}}-1\right)\ge\frac{(n-1)x_1^\gamma
}{\sum_{i=2}^n x_i^\gamma}-1.$$ We may now assume that $x_1=1.$
Hence, we need to show that $$\frac{n}{1+\sum_{i=2}^n
x_i}\left(\frac{1}{G}-1\right)\ge\frac{1}{A}-1,\eqno{(9)}$$ where
$$A=\frac{\sum_{i=2}^nx_i^\gamma}{n-1}\hbox{\ \ and\ \
}G=\ds\prod_{i=2}^n x_i^{\frac{\alpha-1}{n}}.$$ We now choose
$\ds\gamma=\frac{(n-1)(\alpha-1)}{n},$ which implies that $A\ge G$
by the AM-GM Inequality.

For the proof of (9), we will consider two subcases.
\smallskip

{\it Subcase 2.1.} $\ds\prod_{i=2}^n{x_i}\ge 1.$

Then (9) follows from the inequalities $\ds
0\ge\frac{1}{G}-1\ge\frac{1}{A}-1,$ and

$$\frac{\sum_{i=2}^n x_i}{n-1}\ge{\root n-1\of{\ds\prod_{i=2}^n
x_i}}\ge1;\hbox{\ \ \ that is, \  }\frac{n}{1+\sum_{i=2}^n x_i}\le
1.$$

{\it Subcase 2.2.} $\ds\prod_{i=2}^n{x_i}\le 1.$

If $A>1,$ then the left side of (9) is nonnegative, while the
right side is negative, and (9) is thus true. Otherwise,
$\ds\frac{1}{G}-1\ge\frac{1}{A}-1\ge 0,$ and (9) will follow once
we show that $$\frac{n}{1+\sum_{i=2}^n x_i}\ge 1,\hbox{\ \ \ that
is,\ }\frac{\sum_{i=2}^n x_i}{n-1}\le 1.$$ Since
$\ds\alpha\ge2+\frac{1}{n-1}$ (that is $\gamma\ge 1$), the
function $x^\gamma$ (for $x\ge 0$) is convex. Jensen's Inequality
now implies that $$1\ge A \ge\left(\frac{\sum_{i=2}^n
x_i}{n-1}\right)^\gamma,$$ and (9) is proven.
\smallskip

{\bf Remarks.}

a) We already mentioned that the approach in Case 1 is not
applicable in Case 2, since (3) is incorrect when
$\ds\beta\not\in\left[\frac{1}{1-n},1\right].$ Furthermore, it
follows from Subcase 1.1 that, for $\beta>1,$ the opposite
inequality holds. The same inequality follows from Subcase 1.1 for
$\beta\le 1-n$ and $\ds\prod_{i=1}^nx_i=1$ (to see this, consider
the reciprocals of $x_1,\dots,x_n$). On the other side for $n\ge
3,$ $\ds\beta\in\left(1-n,\frac{1}{1-n}\right)$ and
$\ds\prod_{i=1}^n x_i=1$ neither (3) nor its opposite holds true.
To this aim, we note that if $x_1=x^{n-1},$ $\ds
x_2=\dots=x_n=\frac{1}{x},$ then the difference between the left
side and the right side of (3) equals
$$x^{n-1}-x^{\beta(n-1)}+(n-1)\left(\frac{1}{x}-\frac{1}{x^\beta}\right).$$
When $\beta\in(1-n,1)$ and $x\to+\infty,$ this difference tends to
$+\infty,$ while when $\ds\beta<\frac{1}{1-n}$ and $x\to 0^+$ it
tends to $-\infty.$ In particular, when
$\ds\beta\in\left(1-n,\frac{1}{1-n}\right),$ this difference can
have either sign.
\smallskip

b) Let $\ds\gamma=\frac{(n-1)(\alpha-1)}{n}.$ Then (8) holds for
each $\alpha\ge 1$ only when $n=2.$ Indeed, it follows form (9),
which is true for $n=2$ and arbitrary $\alpha\ge 1$ (check!). On
the other side when $\ds\prod_{i=1}^n x_i=1,$ then (8) is
equivalent to (9). When $n\ge 3$, if one chooses
$$x_2=\left((n-\frac{3}{2}\right)^{\frac{1}{\gamma}},\
x_3=\dots=x_n=\left(\frac{1}{2n}\right)^{\frac{1}{\gamma}},$$ then
$G$ and $A$ are less than 1 and independent of $\alpha.$ For
$\alpha\to1^+,$ one has $x_2\to+\infty$ and the left side of (9)
tends to 0, while the right side is a fixed positive number. Thus,
(9) fails for $\alpha$ close to 1, and does (8). For a fixed $n\ge
3,$ one might be interested in finding the least $\alpha_n>1$ that
makes (8) true for each $\alpha\ge\alpha_n$ and to prove that it
fails for $1<\alpha<\alpha_n.$
\smallskip

c) Note that, when $\ds\alpha=2+\frac{1}{n-1},$ we have $\gamma=1$
and (8) follows from (2) by the AM-GM inequality for the numbers
$x_2,\dots,x_n.$
\smallskip

d) Having in mind Proposition 1, one could expect that for each
$\alpha<1$ the opposite inequality to (1) will be satisfied. For
$n=1$ it is trivial, and the reader may easily check it for $n=2.$
Unfortunately, for $n\ge 3$ this is not true. For example, if
$x_1=x^{n-1}>1,$ $\ds x_2=\dots=x_n=\frac{1}{x}$ and
$\alpha\to-\infty,$ the left side of (1) tends to $\ds
n-1-\frac{x^n}{n-1},$ which is positive for $ x^n<(n-1)^2$ ($>1$
for $n>2$). It is interesting to find the least $\alpha_n<1$ such
that the opposite inequality to (1) is true for
$\alpha_n\le\alpha\le 1$ and to prove that it fails for
$\alpha<\alpha_n.$ Here we can consider the following.
\smallskip

{\bf Proposition 2.} {\it Suppose $n\ge 2.$ If
$\frac{1}{1-n}\le\alpha\le 1$ and the positive real numbers
$x_1,x_2,\dots,x_n$ satisfy $\ds\prod_{i=1}^nx_i\ge 1,$ then
$$\sum_{i=1}^n\frac{x_i^\alpha-x_i}{x_1+\dots+x_{i-1}
+x_i^\alpha+x_{i+1}+\dots+x_n}\le 0.$$}

{\it Proof:} This follows from (3) and
$$\frac{x_i^\alpha-x_i}{x_1+\dots+x_{i-1}
+x_i^\alpha+x_{i+1}+\dots+x_n}\le\frac{x_i^\alpha-x_i}{\sum_{i=1}^n
x_i}.$$

{\bf Acknowledgements.} The results of this paper were obtained at
the Summer Research School (2005) of the High School Students'
Institute of Mathematics and Informatics. The author would like to
thank the Institute for the opportunity to take part in the work
of the School. The author is grateful to Professors Oleg
Mushkarov, Peter Popivanov and Ivaylo Kortezov for their useful
remarks and comments on this paper.
\end{document}